\documentstyle[12pt,leqno,amssymb,amsfonts,amscd,graphics,hyperref,tikz]{amsart}  
\hypersetup{hidelinks=true}

\newtheorem{thm}{Theorem}[section]
\newtheorem{lemma}[thm]{Lemma}
\newtheorem{prop}[thm]{Proposition}
\newtheorem{cor}[thm]{Corollary}
\theoremstyle{definition}
\newtheorem{dfn}[thm]{Definition}
\newtheorem{remark}[thm]{Remark}

\begin{document}

\newcommand{\ct}{\cite}
\newcommand{\pr}{\protect\ref}
\newcommand{\su}{\subseteq}
\newcommand{\dn}{\Delta^n}
\newcommand{\pa}{\partial}
\newcommand{\R}{{\Bbb R}}
\newcommand{\Z}{{\Bbb Z}}
\newcommand{\D}{\Delta}
\newcommand{\gt}{\widehat{g}_k}
\newcommand{\gh}{g'}
\newcommand{\p}{{\mathfrak{p}}}
\newcommand{\Y}{\widehat{Y}}
\newcommand{\E}{{\bold{E}}}
\newcommand{\bP}{{\bold{P}}}
\newcommand{\V}{{\mathcal{V}}}
\newcommand{\U}{{\mathcal{U}}}

\newcommand{\k}{{{\ell}}}
\newcommand{\kp}{{\ell'}}
\newcommand{\e}{{\epsilon}}
\newcommand{\tk}{{\E\gt}}

\newcounter{numb}
\numberwithin{equation}{section}

\title[random simplicial complexes in the multi-parameter upper model]{The homology of random simplicial complexes \\ in the multi-parameter upper model}

\author{Michael Farber}
\address{School of Mathematical Sciences, 
	Queen Mary University of London,
	London E1 4NS, UK}
\email{M.Farber@@qmul.ac.uk}
\urladdr{\url{www.qmul.ac.uk/maths/profiles/farberm.html}}

\author{Tahl Nowik}
\address{Department of Mathematics, Bar-Ilan University, 
Ramat-Gan 5290002, Israel}
\email{tahl@@math.biu.ac.il}
\urladdr{\url{www.math.biu.ac.il/~tahl}}
\thanks{M. Farber was partially supported by a grant from the EPSRC}
\date{\today}
\begin{abstract}
We study random simplicial complexes in the multi-parameter upper model. In this model simplices of various dimensions are taken randomly and
 independently, and our random simplicial complex $Y$ is then taken to be the minimal simplicial complex containing this collection of simplices. 

We study the asymptotic behavior of the homology of $Y$ as the number of vertices goes to $\infty$. We observe the following phenomenon asymptotically almost surely. The given probabilities with which the simplices are taken determine a range of dimensions $\k \leq k \leq \kp$ with $\kp \leq 2\k +1$, outside of which the homology of $Y$ vanishes. Within this range, 
the homologies diminish drastically from dimension to dimension. In particular, the homology in the critical dimension $\k$ is significantly the largest.

\end{abstract}

\maketitle

\section{Introduction}\label{intro}

In this work we are interested in random simplicial complexes in the \emph{multi-parameter upper model}. In its most general form it is described as follows. Let $\dn$ be the $n$-dimensional simplex thought of as a simplicial complex. That is, $\dn$ is the set of all non-empty subsets of the set $\{0,\dots,n\}$ of $n+1$ vertices. Let $\Omega_n$ denote the set of all simplicial complexes $Y\su\dn$. Given an assignment of probabilities 
$\{p_\sigma\}_{\sigma\in\dn}$, $0\leq p_\sigma \leq 1$, it induces a probability measure on $\Omega_n$ as follows.
As an intermediate step we randomly select a hypergraph $X \su \dn$, by which we mean an arbitrary subset of $\dn$, not necessarily a simplicial complex.
Each simplex $\sigma \in \dn$ is included in $X$ independently with probability $p_\sigma$,
i.e.\ the probability for obtaining $X$ is 
$\prod_{\sigma \in X} p_\sigma \prod_{\sigma \not\in X} q_\sigma$, where $q_\sigma=1-p_\sigma$.
Now, the simplicial complex produced by this random process is the minimal simplicial complex containing $X$, which we denote by $\overline{X}$. That is, $\overline{X}$ includes all simplices of $X$ and all their faces. 

Though we will not need it in this work, we present an explicit formula for the probability 
$\overline{\bP}(Y)$ that a given simplicial complex $Y$ is obtained in this random process. That is, $\overline{\bP}(Y)$ is the probability that our random hypergraph $X$ satisfies $\overline{X}=Y$.
Let $M(Y)$ denote the set of all maximal simplices in $Y$, then
$$\overline{\bP}(Y)= \prod_{\sigma \not\in Y} q_\sigma\prod_{\sigma \in M(Y)} p_\sigma.$$
This is because $\overline{X}=Y$ iff $X\su Y$ and $X\supseteq M(Y)$. Indeed, if $X\su Y$ then  $\overline{X}\su Y$,
since $Y$ is a simplicial complex, and if $X\supseteq M(Y)$ then $\overline{X}\supseteq Y$, 
since every simplex in $Y$ is contained in a maximal simplex. On the other hand, if $\overline{X}=Y$
then $X\su \overline{X} = Y$ which implies $X\supseteq M(Y)$, since if $\sigma \in M(Y)$ then the only simplex in $Y$ containing $\sigma$ is $\sigma$ itself.

We will be interested in $\{p_\sigma\}$ of a very particular form. We will be given an $(r+1)$-tuple 
$\alpha=(\alpha_0,\dots,\alpha_r) \in [0 ,\infty]^{r+1}$. The probabilities $p_\sigma$ determined by $\alpha\in [0 ,\infty]^{r+1}$ are as follows. Let $j=\dim\sigma$. If $j>r$ or $\alpha_j=\infty$ then $p_\sigma=0$. Otherwise  $p_\sigma=n^{-\alpha_j}$. 

For a fixed $\alpha\in[0 ,\infty]^{r+1}$ we will be interested in the asymptotic behavior of the homologies $\tilde{H}_i(Y;\Z)$ of the random complex $Y$ as $n\to\infty$. In the present work we show that the parameter space $[0 ,\infty]^{r+1}$ may be divided into domains $\U_\k$, where if $\alpha\in \U_\k$ then 
\emph{asymptotically almost surely} (a.a.s.), that is, with probability converging to 1 as $n\to\infty$,
the homology of $Y$ is dominated by that in dimension $\k$. More in detail, there is $\k'\leq 2\k+1$ (depending on $\alpha \in \U_\k$) such that a.a.s.\
the homologies in dimensions $\k \leq k \leq \k'$ diminish drastically from dimension to dimension, and the homologies outside this range vanish.
The dominant dimension $\k$ is named the \emph{critical dimension}.
A specific special case of this phenomenon in the upper model has been studied in \cite{fmn}.

The upper model we have just described is in contrast to the multi-parameter \emph{lower} model, which 
begins with the same random hypergraph $X$ but then produces from it 
the \emph{maximal} simplicial complex \emph{contained} in $X$, which we denote by $\underline{X}$.
That is, a simplex $\sigma$ is in $\underline{X}$ if $\sigma$ and all its faces are in $X$.
The upper and lower models are dual in a clear sense, and so the formula  for 
the probability $\underline{\bP}(Y)$, for obtaining a given simplicial complex $Y$ in the lower model, is dual to that for the upper model. Let $E(Y)$ denote the set of all \emph{minimal} simplices among those \emph{not} in $Y$. An equivalent and geometrically more suggestive definition is 
$E(Y)=\{ \sigma\in\dn : \sigma \not\in Y  \ \hbox{but} \   \pa \sigma \su Y\}$. (This includes the case that $\sigma$ is a vertex $v$ not in $Y$ since then $\pa v = \varnothing \su Y$.)  
We have
$$\underline{\bP}(Y)=\prod_{\sigma \in Y} p_\sigma \prod_{\sigma \in E(Y)} q_\sigma.$$ 
This is because $\underline{X}=Y$ iff $X \supseteq Y$ and $X\cap E(Y) = \varnothing$.
Indeed, if $X \supseteq Y$ then  $\underline{X} \supseteq Y$,
since $Y$ is a simplicial complex, and if $X\cap E(Y) = \varnothing$ then $\underline{X}\su Y$, 
since every simplex not in $Y$ contains a simplex which is minimal among those not in $Y$. On the other hand if $\underline{X}=Y$ then  $X \supseteq\underline{X}= Y$ 
which implies $X\cap E(Y) = \varnothing$, since if $\sigma \in X\cap E(Y)$ then $Y\cup\{\sigma\}$ is a simplicial complex contained in $X$. 

The sets appearing in the formulas for $\overline{P}$ and $\underline{P}$ are $Y,M(Y),E(Y)$, which have clear geometric meaning in terms of $Y$, but make the duality slightly less apparent. From the point of view of duality one may like to add the notation $Y^c$ for the set of all simplices not in $Y$. Then $E(Y)$ is the set of minimal simplices in $Y^c$ and so may be denoted $m(Y^c)$. In terms of this notation we have 
$\overline{\bP}(Y)= \prod_{\sigma \in Y^c} q_\sigma\prod_{\sigma \in M(Y)} p_\sigma$ and
$\underline{\bP}(Y)=\prod_{\sigma \in Y} p_\sigma \prod_{\sigma \in m(Y^c)} q_\sigma$. This makes the duality completely transparent, one formula is obtained from the other by everywhere exchanging 
$Y\leftrightarrow Y^c$, $M\leftrightarrow m$, $p \leftrightarrow q$.
For more on the duality between the upper and lower models see \cite{fmn}. 

The phenomenon of critical dimension also holds in the lower model. This has been established  in \cite{cf}.
The asymptotic behavior observed in the lower model 
resembles that 
of the upper model  but occurs for a different division of the parameter space into domains $\U_\k$.

\section{Definitions and statement of result}\label{statement}

For a fixed integer $r\geq 0$ we study random simplicial complexes of dimension $\leq r$. 
We are given an $(r+1)$-tuple $\alpha=(\alpha_0,\dots,\alpha_r) \in [0 ,\infty]^{r+1}$,
 and an integer $n\geq 0$.
With this data we produce a random hypergraph $X$ by taking each simplex $\sigma$ of dimension
$0\leq i \leq r$ on the vertex set $\{0,\dots,n\}$ independently with probability $n^{-\alpha_i}$ (where
by definition $n^{-\infty}=0$).
Our random simplicial complex $Y$ is then defined to include all the simplices in $X$ and all their faces.

We are interested in the \emph{asymptotic behavior} of our random simplicial complex $Y$, by which we mean the following. We fix the parameters $(\alpha_0,\dots,\alpha_r)$, and we take $n$ to be larger 
and larger. The asymptotic behavior is then described in terms of the following probabilistic notion.

\begin{dfn}\label{aas}
If for every $n$ we have a random object $Z=Z(n)$, and if $T$ is a property that $Z$ may or may not have, then we say $T$ holds \emph{asymptotically almost surely} (a.a.s.) if the probability that $T$ holds converges to 1 as $n\to\infty$.
\end{dfn}

Given $\alpha=(\alpha_0,\dots,\alpha_r) \in [0 ,\infty]^{r+1}$ we define 
$\beta_i=i+1-\alpha_i$ and set 
\begin{equation}\label{beta}
\beta=\beta(\alpha)=\max\{\beta_0,\dots,\beta_r\} \leq r+1.
\end{equation}
We divide our space of parameters $\V=[0 ,\infty]^{r+1}$ into domains $\U_-$ and $\U_0,\dots,\U_r$  where 
$\U_-=\{\alpha\in \V:\ \beta(\alpha)<0\}$ 
and $\U_\k=\{\alpha\in \V:\ \k<\beta(\alpha)<\k+1\}$ for $0\leq \k \leq r$.
The asymptotic behavior of our random simplicial complex $Y$ is as follows.
If $\alpha\in \U_-$, i.e. $\beta<0$, then $Y=\varnothing$ a.a.s., see Proposition \pr{A}.
If $\alpha \in \U_\k$ for $0\leq \k \leq r$, i.e. $\beta>0$ is not an integer and $\k=\lfloor \beta \rfloor$ (the integer part of $\beta$), then 
there is $\k' \leq \lfloor 2\beta \rfloor$  (depending on $\alpha\in\U_\k$) such that a.a.s.
the homologies in dimensions $\k \leq k \leq \k'$ diminish drastically from dimension to dimension, and the homologies outside this range vanish.
The precise details are stated in Theorem \pr{p3}, for which we need one more definition.

\begin{dfn}\label{d}
For two quantities $a=a(n),b=b(n)$, 
\begin{enumerate}
\item We say $a\sim b$ if $\frac{a}{b}\to 1$ as $n\to\infty$.
\item If $a$ is a \emph{random} quantity then we say $a\sim b$ a.a.s. if there is a sequence $\epsilon_n \to 0$ such that
$|\frac{a}{b}-1| <\epsilon_n$ a.a.s., or equivalently, if there is a sequence $c_n$ such that  $|a-b|<c_n$ a.a.s.  and $\frac{c_n}{b}\to 0$.
\end{enumerate}
\end{dfn}

\begin{thm}\label{p3} Let $Y$ be the random simplicial complex in the multi-parameter upper model  determined by  parameters 
$\alpha=(\alpha_0,\dots,\alpha_r)$, and let $\beta=\beta(\alpha)$ be as in (\pr{beta}). Assume $\beta>0$, $\beta\not\in\Z$, and define the critical dimension  $\k=\lfloor \beta \rfloor$. 
We further define the following quantities:
\begin{equation}\label{dd}
d=\sum_{i,\beta_i=\beta}  \frac{1}{(i+1)\k!(i-\k)!}
\end{equation}
\begin{equation}\label{del} 
\nu_k= 2\gamma_k-k \ \ \text{with} \  \
\gamma_k= \max \{\beta_k,\dots,\beta_r\}
\end{equation}
\begin{equation}\label{zero}
\kp=\max\{  k : \ \nu_k \geq 0 \}. 
\end{equation}
Then the following holds a.a.s. (Definition \pr{aas}):
\begin{itemize}
\item $Y$ has full $(\k-1)$-skeleton.
\item $Y$ may be collapsed into its $\kp$-skeleton, \ having  $\k \leq \kp\leq \lfloor2\beta\rfloor\leq 2\k+1$.
\end{itemize}
Let $b_k(Y)$ denote the $k$th Betti number of $Y$, then furthermore, for every sequence $\omega=\omega(n)\to\infty$ the following holds a.a.s.:
\begin{itemize}
\item For $k < \k$, \ $\tilde{H}_k(Y;\Z) = 0$. 
\item For $k=\k$, \  $b_\k(Y) \sim d n^{\beta}$. 
\item For $ \k < k \leq \kp$, \ $b_k(Y) \leq \omega n^{\nu_k}$ \ having  $\nu_{\k+1}<\beta$ and  $\nu_{k+1}\leq\nu_k -1$   for all $k$.
\item For $k > \kp$, \ $\tilde{H}_k(Y;\Z) = 0$. 
\end{itemize}
\end{thm}

We make the following remarks:
\begin{enumerate}
\item The constant $d$ depends on $\alpha$, but it attains only finitely many different values. 
Indeed, $d$  is determined by the (non-empty) set of indices $i$ for which $\beta_i=\beta$.
By Remark \pr{ik} below if $\beta_i=\beta$ then $i\geq \k$.
So we have a stratification of  $\U_\k$ into $2^{r-\k+1}-1$ strata, on each of which $d$ is constant. 
These strata are convex (for this to be meaningful we need to exclude the values $\alpha_i=\infty$). The domain $\U_\k$ itself is connected, but in general it is not convex. 
\item We think of the piecewise linear hypersurfaces $\{\alpha\in \V: \ \beta(\alpha)=\k\}$ that separate between the domains $\U_\k$ in the parameter space $\V$ as \emph{multi-parameter thresholds} for passing from one typical behavior to another. For example, on one side of the piecewise linear hypersurface
$\{\alpha\in \V: \ \beta(\alpha)=\k+1\}$, in $\U_\k$, we have $\tilde{H}_\k(Y;\Z)$ very large, namely 
$b_\k(Y)\sim dn^\beta$ a.a.s. and on the other side of this hypersurface, 
in $\U_{\k+1}$, we have $\tilde{H}_\k(Y;\Z)=0$ a.a.s.
\item We have noted that $\beta\leq r+1$. The boundary case where $\beta=r+1$ is easily understood. In this case $\alpha_r=0$, so every $r$-simplex is included in our random hypergraph $X$ with probability 1, and so $Y$ is the full $r$-skeleton on $\{0,\dots,n\}$  with probability 1.
\item If $\nu_{\k+1}<0$ then  $\k'=\k$, meaning that case (3) of the theorem is empty. That is, all the homology of $Y$ appears only in dimension $\k$. In general we have $\gamma_{\k+1}\leq \beta <\k+1$, whereas the present case $\nu_{\k+1}<0$ means  $\gamma_{\k+1} <\frac{1}{2}(\k+1)$.
\end{enumerate}

The plan of the paper is as follows. In Section \pr{ks} we are interested in the number $f_k$ of $k$-simplices in $Y$ for $k \geq \k$, and give the asymptotic behavior of $f_k$ in Proposition \pr{fk}.
In Section \pr{collapse} we are interested in the Betti numbers $b_k(Y)$ for $k\geq \k$. The asymptotic behavior of these Betti numbers is given in Propositions \pr{p6}, \pr{dec}, \pr{bk}. This is achieved by collapsing $Y$ onto a smaller subcomplex $Y'$.
In Section \pr{homology} we are interested in $\tilde{H}_k(Y;\Z)$ for $k<\k$, showing in Propositions \pr{c3}, \pr{p7} that $\tilde{H}_k(Y;\Z)=0$ a.a.s. This is achieved by a modification of our random model that reduces it to that of Linial-Meshulam. Propositions \pr{p6},  \pr{dec}, \pr{bk}, \pr{c3},  \pr{p7} together constitute  Theorem \pr{p3}.

\section{Counting simplices}\label{ks}

Let $g_i$ denote the number of $i$-simplices in our random hypergraph $X$. Then $g_i$ is a binomial random variable with parameters $\binom{n+1}{i+1}$, $n^{-\alpha_i}$, so 
\begin{equation}\label{egi}
\E g_i=\binom{n+1}{i+1}n^{-\alpha_i}\sim \frac{1}{(i+1)!}n^{i+1-\alpha_i}=\frac{1}{(i+1)!}n^{\beta_i}
\end{equation}
where $\E g_i$ denotes the expectation of $g_i$.
Our first domain $\U_-=\{\beta<0\}$  is easily understood:

\begin{prop}\label{A}
If $\beta<0$ then a.a.s. $X=\varnothing$ and so $Y=\varnothing$.
\end{prop}

\begin{pf}
If $\beta<0$ then $\beta_i < 0$ for all $0\leq i \leq r$.
Markov's inequality gives $\bP\Big(g_i\geq 1\Big)\leq \E(g_i) \to 0$ by (\pr{egi}), i.e. $g_i=0$ a.a.s. for each
$0\leq i \leq r$, so  $X=\varnothing$ a.a.s.
\end{pf}

For the rest of this work we fix an integer $0 \leq \k \leq r$ and an 
$\alpha=(\alpha_0,\dots,\alpha_r)\in \U_\k$.
That is, for our fixed $\k$ we have
\begin{equation}\label{kk}
\k < \beta < \k+1.
\end{equation}
This may also be stated as follows: We assume $0<\beta<r+1$, $\beta\not\in\Z$, and we set 
$\k=\lfloor \beta \rfloor$. For $0 \leq k \leq r$
let $f_k=f_k(Y)$ denote the number of $k$-simplices in $Y$.
Our first goal is  to approximate $f_k$ for $k\geq \k$. Since for every $i \geq k$ each $i$-simplex of $X$ contributes $\binom{i+1}{k+1}$ $k$-simplices to $Y$, 
we have $f_k \leq \sum_{i=k}^r{i+1 \choose k+1}g_i$. 
It is only an inequality since different $i$-simplices may contribute the same $k$-simplex.
This sum will be central in our computations so  we denote 
 $\gt=\sum_{i=k}^r \binom{i+1}{k+1}g_i$ and we have
\begin{equation}\label{0}
f_k \leq \gt.
\end{equation}
By (\pr{egi}) we have
\begin{equation}\label{theta}
\tk=  \sum_{i=k}^r \binom{i+1}{k+1}\binom{n+1}{i+1} n^{-\alpha_i}       
=\binom{n+1}{k+1}\sum_{i=k}^r\binom{n-k}{i-k} n^{-\alpha_i}
\sim \sum_{i=k}^r\frac{n^{\beta_i}}{(k+1)!(i-k)!}.
\end{equation}
Equality holds by the identity $\binom{i+1}{k+1}\binom{n+1}{i+1}       
=\binom{n+1}{k+1}\binom{n-k}{i-k}$ which is true since both sides count the number of pairs of simplices $(\sigma,\tau)$ with $\dim\sigma=k$, $\dim\tau=i$,  $\sigma \su \tau$.

\begin{remark}\label{ik}
If $i$ is such that $\beta_i=\beta$ then $i \geq \k$. 
Indeed, by (\pr{kk}) we have $\k<\beta=\beta_i = i+1-\alpha_i \leq i+1$, so $\k \leq i$.
Note that there may be more than one $i$ such that $\beta_i=\beta$.
\end{remark}

Recall from (\pr{del}) that we define
$\gamma_k=\max\{\beta_k,\dots,\beta_r\}$. 
By Remark \pr{ik} we have 
\begin{equation}\label{gkb}
\gamma_\k=\beta.
\end{equation}
We have $\gamma_k \leq \beta < \k+1 $ by (\pr{kk}), so for  every $k\geq\k$ we have
\begin{equation}\label{k}
\gamma_k < k+1.
\end{equation}

In   Lemma \pr{3}
we give a bound on the difference $\gt-f_k$ for $k \geq \k$.  In the present section it will be used for evaluating $f_k$ via an evaluation of $\gt$. In the next section it will be used for estimating the extent to which the simplices of $X$ overlap.
Note for example that $\gt-f_k=0$ iff every two simplices of $X$ of dimension $\geq k$ intersect in dimension $< k$.

\begin{lemma}\label{3} 
Let $k \geq \k$.
For any sequence $\omega\to\infty$ we have  $\gt-f_k<\omega n^{2\gamma_k-k-1}$  a.a.s.
\end{lemma}

\begin{pf}
For a given $k$-simplex $\sigma$ and $i\geq k$,  there are $\binom{n-k}{i-k}$ $i$-simplices that contain $\sigma$, and so
we have 
$$\bP(\sigma \in Y) = 1-\prod_{i=k}^r(1-n^{-\alpha_i})^\binom{n-k}{i-k}.$$ 
Let  $N_i=\binom{n-k}{i-k}$, $u_i=n^{-\alpha_i}$.
We have 
\begin{multline*}
\bP(\sigma \in Y) = 1-\prod_{i=k}^r(1-u_i)^{N_i}  
=1-\prod_{i=k}^r\Big(1-N_iu_i + \binom{N_i}{2}u_i^2 - \binom{N_i}{3}u_i^3 + \cdots     \Big) \\
\geq 
1-\prod_{i=k}^r\Big(1-N_iu_i + \binom{N_i}{2}u_i^2 \Big)  
\geq 1-\prod_{i=k}^r\Big(1-N_iu_i + (N_iu_i)^2 \Big)
\end{multline*} 
The first inequality holds since for each $i$, $\binom{N_i}{j}u_i^j$ is decreasing in $j$. Indeed, for  $j<N_i$
$$\frac{\binom{N_i}{j+1}u_i^{j+1}}{\binom{N_i}{j}u_i^j} = \frac{N_i-j}{j+1}u_i \leq N_iu_i 
= \binom{n-k}{i-k}n^{-\alpha_i} \leq n^{i-k-\alpha_i}=n^{\beta_i-k-1}\leq n^{\gamma_k-k-1}< 1$$ by (\pr{k}). 
(For the first factor $i=k$ we have $N_k=1$ and the 
inequality for this factor is seen directly.)
So we have 
$\bP(\sigma \in Y)\geq \sum_{i=k}^r N_iu_i -\sum_j T_j$ where each term  $T_j$ is a product of at least two factors of the form $N_iu_i$, and there are less than $3^{r+1}$ such terms. 
Again using  $N_iu_i \leq n^{\gamma_k-k-1}< 1$ we get 
$\bP(\sigma \in Y) \geq 
\Big(\sum_{i=k}^r N_iu_i\Big)\Big(1-cn^{\gamma_k-k-1}\Big)$
for a constant $c>0$.
Thus   $$\E f_k =\binom{n+1}{k+1}  \bP(\sigma \in Y) \geq \binom{n+1}{k+1}\Big(\sum_{i=k}^rN_iu_i\Big)\Big(1-cn^{\gamma_k-k-1}\Big) = \tk(1-cn^{\gamma_k-k-1})$$
by  (\pr{theta}), which may be rewritten as $\E (\gt-f_k) \leq c n^{\gamma_k-k-1}  \tk$.
Now, by (\pr{theta}) we have 
$\tk\leq c'n^{\gamma_k}$ for some $c'>0$, so together
$\E (\gt-f_k) \leq cc' n^{2\gamma_k-k-1}$.
Since $\gt-f_k \geq 0$ we may use Markov's inequality 
$$\bP\Big(\gt-f_k \geq \omega n^{2\gamma_k-k-1}\Big) 
\leq \frac{\E (\gt-f_k)}{\omega n^{2\gamma_k-k-1}} 
  \to 0.$$
\end{pf}

We will have two occasions to use the following lemma, with different choices of coefficients.

\begin{lemma}\label{cheb}
Given $a_k,\dots,a_r$ with $a_i > 0$, let $g^a=\sum_{i=k}^r a_ig_i$. We have:
\begin{enumerate}
\item If $\gamma_k >0$ then $g^a \sim A^a n^{\gamma_k}$ a.a.s. with 
$A^a=\sum_{i\geq k,\beta_i=\gamma_k} \frac{a_i}{(i+1)!}$
\item If $\gamma_k=0$ then for every sequence $\omega\to\infty$ we have $g^a \leq \omega$ a.a.s.
\item If $\gamma_k <0$ then $g^a=0$ a.a.s.
\end{enumerate}
\end{lemma}

\begin{pf} We prove in opposite order: 
(3) As in Proposition \pr{A}. 

(2) If there are $k \leq i \leq r$ with $\beta_i<0$ then (3) applies to them. For $i$ with $\beta_i=0$, i.e. $\alpha_i=i+1$,  $g_i$ is a binomial random variable with parameters $\binom{n+1}{i+1}$, $n^{-(i+1)}$. We have $\binom{n+1}{i+1} n^{-(i+1)} \to \frac{1}{(i+1)!}$, so the
distribution of $g_i$ converges to a Poisson distribution, the claim follows. 

(1) If there are $k \leq i \leq r$ with $\beta_i\leq 0$ then (3) or (2) apply to them, taking $\omega=n^\e$ with  $0<\e<\gamma_k$ when (2) applies.
For $i$ with $\beta_i>0$, $g_i$ is a binomial random variable with parameters $\binom{n+1}{i+1}$, $n^{-\alpha_i}$, so by Chebyshev's inequality  we have
$$\bP\Big(|g_i - \E g_i|  \geq n^{\frac{2}{3}\beta_i}  \Big) 
\leq \frac{{\bold{Var}}(g_i)}{n^{\frac{4}{3}\beta_i}}
= \frac{\binom{n+1}{i+1}n^{-\alpha_i}(1-n^{-\alpha_i})}{n^{\frac{4}{3}\beta_i}} 
\leq \frac{(n+1)^{\beta_i}}{n^{\frac{4}{3}\beta_i}}
\to 0.$$
By (\pr{egi}) we have $\E g_i \sim \frac{1}{(i+1)!}n^{\beta_i}$ so $\frac{n^{\frac{2}{3}\beta_i}}{\E g_i} \to 0$  so $g_i \sim \E g_i$ a.a.s.
 so  $g_i \sim \frac{1}{(i+1)!}n^{\beta_i}$ a.a.s.
The claim follows.
\end{pf}

We arrive at the main result of this section.

\begin{prop}\label{fk} Let $f_k$ be the number of $k$-simplices in $Y$.
For $k \geq \k$ we have:
\begin{enumerate}
\item If $\gamma_k >0$ then $f_k \sim D_k n^{\gamma_k}$ a.a.s. with 
$D_k=\sum_{i\geq k,\beta_i=\gamma_k} \frac{1}{(k+1)!(i-k)!}$
\item If $\gamma_k=0$ then for every sequence $\omega\to\infty$ we have $f_k \leq \omega$ a.a.s.
\item If $\gamma_k <0$ then $f_k=0$ a.a.s.
\end{enumerate}
\end{prop}

\begin{pf}
Take $a_i = \binom{i+1}{k+1}$ in Lemma \pr{cheb}, giving $g^a=\gt$.

(1) We have $A^a=D_k$ so by Lemma \pr{cheb}(1) we have 
$\gt \sim D_k n^{\gamma_k}$ a.a.s. 
Now let $\omega = n^{\frac{1}{2}(k+1-\gamma_k)}$ then $\omega\to\infty$ 
 by (\pr{k}), and  
 $\gt-f_k<\omega n^{2\gamma_k-k-1} = \frac{1}{\omega}n^{\gamma_k}$ a.a.s.\   by Lemma \pr{3}. This gives $f_k \sim D_k n^{\gamma_k}$ a.a.s.

For (2),(3) use $f_k \leq \gt$ and Lemma \pr{cheb}(2),(3).
\end{pf}
We mention that $f_k$ for $k<\k$ is given by Proposition \pr{p7} below (which is part of Theorem \pr{p3}).
Namely, for $k<\k$ we have $f_k = \binom{n+1}{k+1}$ a.a.s.

\section{Collapsing simplices}\label{collapse}

Recall $X$ is the random hypergraph that produces our random simplicial complex $Y$, and let $k \geq \k$.
Let $Y_k = \{ \sigma\in Y : \  \dim\sigma=k\}$.
Let $X_{k^+} = \{\tau\in X : \ \dim\tau \geq k \}$, which are
all the simplices in $X$ that contribute $k$-simplices to $Y$.

\begin{dfn}\label{kgood} Let $\tau\in X$ and  $k\geq\k$.
We  say that $\tau$  is \emph{$k$-good} if $\tau\in X_{k^+}$ and for any other  $\tau'\in X_{k^+}$
we have $\dim (\tau \cap \tau') < k$. We say that $\tau$ is \emph{$k$-bad} if $\tau\in X_{k^+}$ and $\tau$ is not $k$-good. 
\end{dfn}

To avoid confusion we emphasize that being $k$-good or $k$-bad is a property of simplices in $X$, not in $Y$. We also note that if $\k \leq k \leq i$ and $\tau$ is a $k$-good $i$-simplex, then $\tau$ is also $j$-good for every $k \leq j \leq i$.

\begin{dfn}\label{good}
Let $\tau\in X$.
\begin{enumerate}
\item We  say $\tau$ is \emph{good} if $\tau$ is $k$-good for some $k\geq\k$.
\item If $\tau$ is good then we denote by $G(\tau)$ the minimal $k\geq\k$ such that $\tau$ is $k$-good. 
\end{enumerate}
\end{dfn}

\begin{lemma}\label{col}
The simplicial complex $Y$ may be collapsed onto a subcomplex $Y' \su Y$
such that for every good simplex $\tau$, 
if $i=\dim\tau$ and $k=G(\tau)$ then:
\begin{enumerate}
\item All $j$-faces of $\tau$ with $j>k$ are removed.
\item All $j$-faces of $\tau$ with $j<k$ remain.
\item Precisely $\binom{i}{k}$ of the $k$-faces of $\tau$ remain. 
\end{enumerate}
\end{lemma}

\begin{pf} We describe the collapse corresponding to each good simplex. 
Let $\tau=\{v_0,\dots,v_i\}$ be a good $i$-simplex with $G(\tau)=k$. 
By definition $i \geq k$ and assume first that $k>0$.
Let $\hat{\tau}$ denote the subcomplex of $Y$ consisting of $\tau$ and all its faces, and let
 $\tau_0=\{v_1,\dots,v_i\}$ be the $(i-1)$-face of $\tau$ opposite to $v_0$. For
$j\leq i-1$ let $\Delta_j$ denote the $j$-skeleton of $\hat{\tau}_0$, and let $v_0 * \Delta_j$ denote 
the cone over $\Delta_j$ with vertex $v_0$.
We collapse $\hat{\tau}$ onto $v_0 * \Delta_{k-1}$ doing it step by step 
$$\hat{\tau}=v_0 * \Delta_{i-1} \ \longrightarrow \ v_0 * \Delta_{i-2} \ \longrightarrow \ \cdots \ \longrightarrow \ v_0 * \Delta_{k-1}.$$ 
For the collapse $v_0 * \Delta_j\to v_0 * \Delta_{j-1}$ we go over all $j$-simplices  $\rho \in \Delta_j$, and for each such $\rho$ we remove the pair of simplices $\rho,  \{v_0\} \cup \rho$. See Figure \pr{f} where  $i=3$, $k=1$.
For all this to be a collapse, the following needs to hold at each  stage.
If $\hat{\tau}$ has already been collapsed onto $v_0 * \Delta_j$ with $j\geq k$, and
$\rho \in \Delta_j$ is a $j$-simplex, then
$\{v_0\} \cup \rho$ is the only simplex that strictly contains $\rho$. 
This is indeed true since $\rho$ is a maximal simplex in 
$\Delta_j$, and since $\tau$ was a $k$-good simplex and $\dim\rho=j\geq k$.

We have collapsed 
 $\hat{\tau}$ onto $v_0 * \Delta_{k-1}$. We note that $\Delta_{k-1}$  includes $\binom{i}{k}$ $(k-1)$-simplices, and so 
$v_0 * \Delta_{k-1}$  includes $\binom{i}{k}$ $k$-simplices as claimed in (3), completing the case $k>0$.
If $k=0$ then  $\tau$ is disjoint from all other simplices of $X$, so $\hat{\tau}$ may be collapsed to $\{v_0\}$  and the claim holds as well.
\end{pf}

\begin{figure}
[h]
\begin{center}
\begin{tikzpicture}[scale=0.6]

\begin{scope}

\fill[gray!60]
(0,0)--(3.5,1.5)--(0.5,4)--(0,0)
(3.5,1.5)--(0.5,4)--(4,4.5)--(3.5,1.5)
; 
\draw (0,0) -- (3.5,1.5) 
(0,0) -- (0.5,4)
(3.5,1.5) -- (0.5,4) 
(0.5,4) --  (4,4.5)
(4,4.5)  -- (3.5,1.5)
;
\filldraw (0,0)  circle(3pt) 
(3.5,1.5) circle(3pt)
(0.5,4)  circle(3pt) 
(4,4.5)  circle(3pt) 
;
\draw[->][line width=1.5] (5,2.5)--(7,2.5);
\end{scope}

\begin{scope}[shift={(8,0)}]

\fill[gray!30]
(0,0)--(3.5,1.5)--(4,4.5)--(0,0)
(0,0)--(0.5,4)--(4,4.5)--(0,0)
; 

\draw (0,0) -- (4,4.5)
;

\fill[gray!60]

(0,0)--(3.5,1.5)--(0.5,4)--(0,0)
; 

\draw (0,0) -- (3.5,1.5) 
(0,0) -- (0.5,4)
(3.5,1.5) -- (0.5,4) 
(0.5,4) --  (4,4.5)
(4,4.5)  -- (3.5,1.5)
;

\filldraw (0,0)  circle(3pt) 
(3.5,1.5) circle(3pt)
(0.5,4)  circle(3pt) 
(4,4.5)  circle(3pt) 
;

\draw[->][line width=1.5] (5,2.5)--(7,2.5);

\end{scope}

\begin{scope}[shift={(16,0)}]

\draw (0,0) -- (4,4.5)
 (0,0) -- (3.5,1.5) 
(0,0) -- (0.5,4)
;
\filldraw (0,0)  circle(3pt) 
(3.5,1.5) circle(3pt)
(0.5,4)  circle(3pt) 
(4,4.5)  circle(3pt) 
;
\end{scope}

\end{tikzpicture}
\end{center}
\caption{$\hat{\tau}=v_0*\Delta_2 \ \longrightarrow \ v_0*\Delta_1 \ \longrightarrow \ v_0*\Delta_0$}
\label{f}
\end{figure}
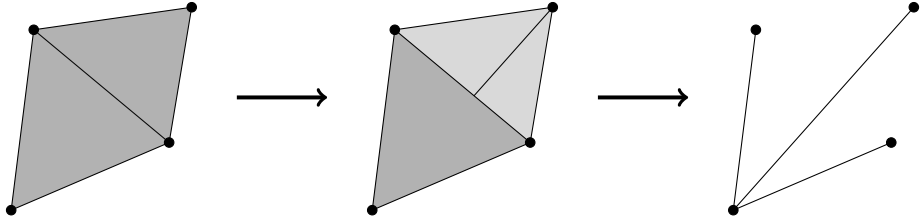

In view of Lemma \pr{col}, for each $k \geq \k$  we would like to have a bound on the number of $k$-bad simplices.

\begin{dfn}
Let $B_k=B_k(X)$ denote the number of $k$-bad simplices.
\end{dfn}

\begin{lemma}\label{p2}
Let $k \geq \k$.
For any sequence $\omega\to\infty$ we have  $B_k \leq \omega n^{2\gamma_k-k-1}$ a.a.s. 
\end{lemma}

\begin{pf}
We show $B_k\leq 2(\gt-f_k)$ which together with Lemma \pr{3} establishes our claim.
For $\sigma\in Y_k$, let $M_\sigma=\{\tau\in X_{k^+} :  \ \tau\supseteq\sigma   \}$ and let
$m_\sigma=|M_\sigma|$ (the number of elements in $M_\sigma$), so $m_\sigma \geq 1$.
We claim $\gt = \sum_{\sigma\in Y_k} m_\sigma$. 
Indeed, both sides of the equality count the number of pairs 
$(\sigma, \tau) \in Y_k \times X_{k^+}$ with $\sigma \su \tau$.
On the other hand $f_k = \sum_{\sigma\in Y_k} 1$. So 
$$\gt-f_k= \sum_{\sigma\in Y_k}(m_\sigma - 1) = \sum_{\sigma\in Y_k, m_\sigma \geq 2}(m_\sigma - 1).$$
Now if $\tau\in X_{k^+}$ is a $k$-bad simplex and $\tau'\in X_{k^+}$, $\tau'\neq\tau$,  is such that  $\dim (\tau \cap \tau') \geq k$ then $\tau \cap \tau'$ contains a $k$-simplex $\sigma$. Thus, a simplex 
$\tau\in X_{k^+}$ is $k$-bad iff $\tau$ contains a $k$-simplex $\sigma$ with  $m_\sigma \geq 2$, that is, the set of all $k$-bad simplices is $\bigcup_{\sigma\in Y_k,  m_\sigma \geq 2} M_\sigma$.
So we have
$$B_k = \Big|\bigcup_{\sigma\in Y_k,   m_\sigma \geq 2} M_\sigma\Big| \leq \sum_{\sigma\in Y_k,  m_\sigma \geq 2} m_\sigma  \leq \sum_{\sigma\in Y_k, m_\sigma \geq 2}2(m_\sigma - 1)=2 (\gt-f_k).$$
\end{pf}

Recall $b_k(Y)$ denotes the $k$th Betti number of $Y$.

\begin{prop}\label{p6} 
For $k>\k$
let $\nu_k$ be as in (\pr{del}), then the following holds:
\begin{enumerate}
\item Given any sequence $\omega\to\infty$ we have
  $b_k(Y) \leq \omega n^{\nu_k}$ a.a.s.
\item If $\nu_k <0$ then $Y$ is collapsible into its $(k-1)$-skeleton a.a.s. so $\tilde{H}_k(Y;\Z)=0$ a.a.s.
\end{enumerate}
\end{prop}

\begin{pf} 
Let $Y'$ denote the collapsed complex given by Lemma \pr{col}
and let $f'_k$ denote the number of $k$-simplices in $Y'$.

(1) Let $\sigma \in Y' \su Y$ be a $k$-simplex, then there is $\tau \in X_{k^+} \su  X_{(k-1)^+}$ such that $\tau \supseteq \sigma$, and we claim that  $\tau$ is $(k-1)$-bad. Indeed, otherwise $\tau$ is $(k-1)$-good so $G(\tau)\leq k-1$ and so 
by Lemma \pr{col}(1) $\sigma$ has been removed during the collapse.
This gives $f'_k \leq \binom{r+1}{k+1}B_{k-1}$ since an  $i$-simplex contains $\binom{i+1}{k+1}$ $k$-simplices and $\binom{i+1}{k+1}\leq \binom{r+1}{k+1}$. 
By Lemma \pr{p2} we 
get $f'_k  \leq \omega n^{2\gamma_{k-1}-k}$ a.a.s.

To obtain (1) we need to replace $\gamma_{k-1}$ with $\gamma_k$ in the last inequality. To achieve this we also look at 
the hypergraph $\widetilde{X}$
 obtained from $X$ by deleting all $(k-1)$-simplices. Let $\widetilde{Y},\widetilde{Y}',\widetilde{f}'_k$ etc.\ be the corresponding objects. 
The possibilities for collapsing simplices of dimension $\geq k$ in $Y$ and $\widetilde{Y}$ are identical, so we look at the collapse assigned to $\widetilde{Y}$ by Lemma \pr{col} and apply it to the simplices of dimension
$\geq k$ in $Y$. (Some $(k-1)$-simplices of $Y$ are also  removed in this process.)
We thus obtain a collapse of $Y$ after which $\widetilde{f}'_k$ $k$-simplices remain, so $b_k(Y)\leq\widetilde{f}'_k$.
The random model that starts with $X$ and then deletes all $(k-1)$-simplices to obtain $\widetilde{X}$ is equivalent to our usual model only with a different assigned probability in dimension $k-1$, namely,  
 $n^{-\alpha_{k-1}}$ is replaced with 0, i.e.\
$\widetilde{\alpha}_{k-1}=\infty$. This gives $\widetilde{\beta}_{k-1}=-\infty$ so $\widetilde{\gamma}_{k-1}=\gamma_k$.
Being equivalent to the usual model,  the bound concluding the previous paragraph applies, and we get
$b_k(Y)\leq\widetilde{f}'_k  \leq \omega n^{2\widetilde{\gamma}_{k-1}-k}=\omega n^{2\gamma_k-k}=\omega n^{\nu_k}$ a.a.s.

(2) We continue  looking at the collapse of $Y$  induced by that of $\widetilde{Y}$.
We have $\widetilde{f}'_k  \leq \omega n^{\nu_k}$ a.a.s. and if $\nu_k <0$ then we can take $\omega=n^\e$
with $\e>0$ such that $\nu_k + \e<0$, so  $\omega n^{\nu_k} \to 0$.
But $\widetilde{f}'_k$ is a sequence of  non-negative integers 
so in fact $\widetilde{f}'_k=0$ a.a.s.
\end{pf}

We remark about the proof above, that the difference between  $f'_k$ and $\widetilde{f}'_k$ is only due to our specific definition of $(k-1)$-good, which in turn determines the 
specific collapse of Lemma \pr{col}. 
There may be an $i$-simplex $\tau$ with $i\geq k$ which is $(k-1)$-good in $\widetilde{X}$, but there is a $(k-1)$-simplex $\sigma$ in $X$ with $\sigma \su \tau$ so $\tau$ is $(k-1)$-bad in $X$.
In our modified collapse of $Y$ using $\widetilde{X}$,  such $\tau$ gets to be collapsed.
The collapse of Lemma \pr{col} with no modification will be used in the proof of Proposition \pr{bk}, followed by a discussion of its efficiency.

We now show that the exponents $\nu_k$  and the dimension $\k'$  satisfy 
the properties stated in 
Theorem \pr{p3}.

\begin{prop}\label{dec}
The exponents $\nu_k$ defined in (\pr{del}) and $\k'$ defined in (\pr{zero}) satisfy the following:
\begin{enumerate}
\item $\nu_{\k+1} < \beta$.
\item $\nu_{k+1} \leq \nu_k - 1$ for all $k$.
\item $\k \leq \kp \leq \lfloor 2\beta \rfloor$.
\end{enumerate}
\end{prop}

\begin{pf}

(1) $\nu_{\k+1} = 2\gamma_{\k+1}-\k-1 \leq 2\beta-\k-1< \beta$ 
by (\pr{kk}). 

(2) $\nu_{k+1}=2\gamma_{k+1}-k-1\leq 2\gamma_k-k-1 = \nu_k-1$.

(3) We have $\nu_\k = 2\gamma_\k - \k=2\beta - \k > \beta >0$ by  (\pr{gkb})
and (\pr{kk}),
so $\k \leq \kp \leq r$. (We use $\nu_\k$  only here, otherwise only $\nu_k$ with $k > \k$ is of interest.)
Now assume $\lfloor 2\beta \rfloor <r$, otherwise we are done. If $\k<\beta<\k+\frac{1}{2}$ then $\lfloor 2\beta \rfloor = 2\k$
and  $\nu_{\k+1} \leq    2\beta-\k-1   <\k$, so by iterating  (2) $\k$ times we
have $\nu_{2\k+1}<0$ so $\kp\leq 2\k=\lfloor 2\beta \rfloor$.
($\nu_{2\k+1}$ is indeed defined since $2\k=\lfloor 2\beta \rfloor<r$.)
Similarly,
if $\k+\frac{1}{2}\leq\beta<\k+1$ then $\lfloor 2\beta \rfloor = 2\k+1$
and   $\nu_{\k+1}\leq  2\beta-\k-1  <\k+1$ so by iterating (2) $\k+1$ times 
we have $\nu_{2\k+2}<0$ so $\kp\leq 2\k+1=\lfloor 2\beta \rfloor$. ($\nu_{2\k+2}$ is indeed defined since $2\k+1=\lfloor 2\beta \rfloor<r$.)
\end{pf}

We now  evaluate $b_\k(Y)$, again via the collapsed complex $Y'$. By Lemma \pr{col}(3)
every $\k$-good $i$-simplex $\tau$ contributes $\binom{i}{\k}$ $\k$-simplices after being collapsed, as opposed to the $\binom{i+1}{\k+1}$ $\k$-simplices contained in $\tau$ before the collapse. We would thus like to have  a ``collapsed version'' of our quantity $\widehat{g}_\k$ where we replace the coefficients $\binom{i+1}{\k+1}$ by  $\binom{i}{\k}$.

\begin{lemma}\label{cv}
Let $\gh=\sum_{i=\k}^r \binom{i}{\k}g_i$, 
then 
 $\gh \sim dn^\beta$ a.a.s. with $d$ given in (\pr{dd}).
\end{lemma}

\begin{pf}
Take $k=\k$ and $a_i = \binom{i}{\k}$ in Lemma \pr{cheb}, then $g^a=\gh$. 
By (\pr{gkb}) we have $\gamma_\k=\beta>0$, so case (1) of Lemma \pr{cheb} applies.
Finally, by Remark \pr{ik}, $i\geq\k$ for every $i$ such that $\beta_i=\beta$, so $A^a=d$.

\end{pf}

\begin{prop}\label{bk}
 $b_\k(Y)\sim dn^\beta$ a.a.s.
with $d$ given in (\pr{dd}).
\end{prop}

\begin{pf} 
As before, let $Y'$ denote the collapsed complex given by Lemma \pr{col}
and let $f'_k$ denote the number of $k$-simplices in $Y'$.
For $i \geq \k$ denote by $g^G_i$ the number of $\k$-good $i$-simplices.
By definition of $G(\tau)$, if $\tau$ is $\k$-good  then $G(\tau)=\k$, and so 
by Lemma \pr{col}(3) every $\k$-good $i$-simplex contributes $\binom{i}{\k}$ $\k$-simplices to $Y'$.   By definition of $\k$-good simplices there is no overlap in these contributions, so we have  
\begin{equation}\label{g}
\sum_{i=\k}^r \binom{i}{\k}g^G_i \leq f'_\k \leq \sum_{i=\k}^r \binom{i}{\k}g^G_i+\binom{r+1}{\k+1}B_\k.
\end{equation}
We further note that $g_i - B_\k \leq g^G_i \leq g_i$ for every $i$. Substituting this into (\pr{g}) gives 
$\gh-cB_\k \leq f'_\k \leq \gh+cB_\k$ for some  $c>0$, where
by Lemma \pr{cv} we have $\gh \sim  dn^\beta$ a.a.s. 
By (\pr{kk}) we have  $2\beta - \k -1 <\beta$ and take $\e>0$ so that
 $\delta=2\beta - \k -1+\e < \beta$.   
Taking $\omega = n^\e$ in Lemma \pr{p2} we get $B_\k\leq n^\delta$ a.a.s. 
since $\gamma_\k=\beta$ by (\pr{gkb}).
Together we get that $f'_\k \sim dn^\beta$ a.a.s.

Now, in the proof of Proposition \pr{p6} we noted  that  $f'_{\k+1}  \leq cB_\k$ for some  $c>0$ so $f'_{\k+1}\leq c n^\delta$ a.a.s. 
Furthermore $f'_{\k-1} \leq \binom{n+1}{\k} \leq (n+1)^\k$ 
and $\k < \beta$ by (\pr{kk}).  So using 
$f'_\k-f'_{\k+1}-f'_{\k-1} \leq b_\k(Y) \leq f'_\k$ we get $b_\k(Y) \sim d n^\beta$ a.a.s. 
\end{pf}

Our collapse pattern of Lemma \pr{col} involves certain choices that may seem arbitrary and perhaps not as efficient as possible. We can now see that in dimension $\k$ only negligible further collapse may be possible. Indeed, by the proof of Proposition \pr{bk}
 the number of $\k$-simplices in our particular collapse satisfies $f'_\k \sim dn^\beta$ a.a.s. and let
$f_\k''(Y)$ denote the minimal  number of $\k$-simplices in any collapse of  $Y$. Then 
$b_\k(Y) \leq f_\k''(Y)\leq f'_\k$ and so  also $f_\k''(Y)\sim d n^\beta$ a.a.s.

We would like to compare the behavior of $b_k(Y)$ described in Propositions \pr{p6}, \pr{bk}
 to that of $f_k(Y)$ described in Proposition \pr{fk}. For $k=\k$ we have the same exponent $\beta$ by (\pr{gkb}). As to the coefficient, in general $d<D_\k$ since 
each term in the sum for  $d$  is $\frac{\k+1}{i+1}$ times the corresponding term in the sum for  $D_\k$. This reflects the fact that when collapsing an $\k$-good $i$-simplex $\tau$, a fraction 
$\frac{\binom{i}{\k}}{\binom{i+1}{\k+1}}=\frac{\k+1}{i+1}$ of the $\k$-faces of $\tau$ survive the collapse.

For $k>\k$, the difference $e_k=\nu_k-\gamma_k$ between the corresponding exponents is 
$e_k=\nu_k-\gamma_k = \gamma_k-k$ which is negative, and drastically more so from dimension to dimension.
Indeed, as in the proof of Proposition \pr{dec}(1),(2) we get $e_{\k+1}<0$ and $e_{k+1}\leq e_k -1$
 for all $k$. This reflects the increasing proportion of collapse that takes place as we go up in the dimensions.

This completes our analysis regarding the homologies of $Y$ in dimensions $k \geq \k$. 
The homologies $\tilde{H}_k(Y;\Z)$ for $k<\k$ are addressed in the next section.
 At this point the collapsed complex $Y'$  has completed its role in our computations and we return to our original random complex $Y$.

\section{The homology $\tilde{H}_k(Y;\Z)$ for $k<\k$}\label{homology}

In case $\k=0$ our analysis is already complete, so we assume $\k>0$.
We analyze $\tilde{H}_{\k-1}(Y;\Z)$ by reduction to the $\k$-dimensional Linial-Meshulam model appearing in \cite{mw}. We start with the full $(\k-1)$-skeleton $K$ on the vertex set $\{0,\dots,n\}$ and use our random hypergraph $X$ to add $\k$-simplices to $K$ by a certain rule presented below. This modified model for a random complex produces an $\k$-complex that we denote $\Y$. 
We will make sure that the $\k$-simplices are added independently with probabilities bounded below by 
$cn^a$ with $a>-1$.
It then follows from Theorem 1.1 of \cite{hkp} that $\tilde{H}_{\k-1}(\Y;\Z)=0$ a.a.s.
We will use this to deduce our desired results regarding our original random complex $Y$.

For this construction,
choose one index $i$ such that $\beta_i=\beta$ and fix it for the rest of this section. We have $i\geq \k$ by 
Remark \pr{ik}, and let $X_i = \{\tau \in X : \ \dim\tau = i\}$. 
In our modified random model we use  $X_i$ for adding $\k$-simplices to $K$. 
But note that if we add to $K$ all $\k$-faces of the simplices in $X_i$ then, if $i>\k$,  the  $\k$-simplices will not be added independently (since for example we would have $\bP(f_\k=1)=0$
and $\bP\big(f_\k=\binom{i+1}{\k+1}\big)>0$).
To circumvent this problem, we will add only one $\k$-face of each $\tau\in X_i$. 
In order that the $\k$-faces will be added with sufficiently large probability,
we wish to have a function that chooses an $\k$-face from each $i$-simplex in a way that 
every $\k$-simplex is chosen by sufficiently many $i$-simplices. 
Let $S_j=S_j(n)$ denote the set of all $j$-simplices on our set $\{0,\dots,n\}$ of vertices.

\begin{lemma}\label{15}

For every sufficiently large $n$ there exists a function  
$h:S_i \to S_\k$  satisfying the following two properties:
\begin{enumerate}
\item $h(\tau)\su\tau$ for every $\tau\in S_i$
\item $|h^{-1}(\sigma)|\geq \frac{\binom{n-\k}{i-\k}}{2\binom{i+1}{\k+1}}$
for every  $\sigma\in S_\k$
\end{enumerate}
\end{lemma}

\begin{pf} Assume first that $i>\k$.
We prove existence of a function $h$ satisfying (1) and (2) using the probabilistic method. 
For each $\tau \in  S_i$ we choose $h(\tau)$  randomly from among the $\binom{i+1}{\k+1}$ $\k$-faces of $\tau$, with equal probabilities and independently. The function $h$ satisfies property (1) by definition.
If we show that there is a  positive probability that $h$ satisfies property (2), then
there must exist at least one such function $h$.

Fix one $\sigma\in S_\k$. Each $\tau \in S_i$ that contains $\sigma$ will choose $\sigma$ to be $h(\tau)$
with probability $\frac{1}{\binom{i+1}{\k+1}}$, independently. 
So $F=|h^{-1}(\sigma)|$ is a binomial random variable with parameters 
$\binom{n-\k}{i-\k}, \frac{1}{\binom{i+1}{\k+1}}$, and we have 
$\E F=\frac{\binom{n-\k}{i-\k}}{\binom{i+1}{\k+1}}$. 
By Chernoff's bound (see e.g. Theorem 2.1 of \cite{jlr}),
$\bP\Big(F< \E F-R\Big)\leq \exp\Big(-\frac{R^2}{2\E F} \Big)$. Taking $R=\frac{1}{2}\E F$ we get
$\bP(F<\frac{1}{2}\E F)\leq \exp(-\frac{1}{8}\E F) \leq  \exp(-cn^{i-\k})$.
This is true for every $\sigma\in S_\k$, thus the probability that there exists \emph{some} $\sigma\in S_\k$ with
$|h^{-1}(\sigma)|<\frac{\binom{n-\k}{i-\k}}{2\binom{i+1}{\k+1}}$ is at most
$\binom{n+1}{\k+1}\exp(-cn^{i-\k})$. We assumed here that $i>\k$, so for sufficiently large $n$ this probability is strictly less than 1, and so for each such $n$ there must exist a function $h$ with the desired property.

In case  $i=\k$  take $h$ to be the identity, giving  $|h^{-1}(\sigma)|=1\geq\frac{1}{2}$.
\end{pf}

For each sufficiently large $n$ we choose one function $h_n$ provided by Lemma \pr{15}. We will use this sequence of functions $h_n$ to define our modified random model. To avoid confusion we emphasize that the probabilistic argument in the proof of Lemma \pr{15} was only a method for proving that functions $h_n$ with the desired properties exit. But once the sequence of functions $h_n$ is chosen, they are fixed once and for all
and are \emph{not} random objects in our  modified random model. 
Accordingly, the sets $h_n^{-1}(\sigma)$ are fixed beforehand once and for all.

Finally, our modified random complex $\Y$ is defined as follows. 
Recall $X_i = \{\tau \in X : \ \dim\tau = i\}$ where  $X$ is the random hypergraph that produces our random simplicial complex $Y$.
We  start with the full $(\k-1)$-skeleton $K$ on the vertex set $\{0,\dots,n\}$, and  for each $\tau\in X_i$ we add  to $K$
 the $\k$-simplex $h_n(\tau)$.

\begin{lemma}\label{p8}
$\tilde{H}_{\k-1}(\Y;\Z)=0$ a.a.s.
\end{lemma}

\begin{pf}
An $\k$-simplex is included in $\Y$ iff one of the $i$-simplices in $h_n^{-1}(\sigma)$ is chosen in the random process defining $X$. Since the sets $h_n^{-1}(\sigma)$ are disjoint, it follows that the $\k$-simplices are included in $\Y$  independently.
We now evaluate the probability that a given $\k$-simplex $\sigma$ is included in $\Y$.
Denote $N=\binom{n-\k}{i-\k}$, $u= n^{-\alpha_i}$, and 
$N_\sigma= |h_n^{-1}(\sigma)|$, then $eN \leq N_\sigma \leq N$ 
with $e=\frac{1}{2\binom{i+1}{\k+1}}$, by Lemma \pr{15}.
We have 
$$\bP(\sigma\in \Y) = 1-(1-u)^{N_\sigma} \geq N_\sigma u -(N_\sigma u)^2$$
since the terms in the alternating binomial sum are decreasing as in  
the proof of Lemma \pr{3}, using $N_\sigma u \leq Nu \leq n^{i-\k-\alpha_i} \leq n^{\beta-\k-1}<1$.
(Here again, in case $i=\k$ we have $N_\sigma =1$ and the inequality is  seen directly.)
Thus we have for sufficiently large $n$
$$\bP(\sigma\in \Y) \geq eNu - (Nu)^2 = (e-Nu)Nu \geq cNu \geq c' n^{i-\k-\alpha_i}= c' n^{\beta-\k-1}$$ 
since $i$ was chosen such that $\beta_i=\beta$. By  (\pr{kk}) we have  $\beta -\k -1 > -1$, so  as mentioned in the opening paragraph of this section, 
it follows from Theorem 1.1 of \cite{hkp} that $\tilde{H}_{\k-1}(\Y;\Z)=0$ a.a.s.
\end{pf}

Returning to our original random complex $Y$ we get the following.

\begin{cor}\label{c2}
$\tilde{H}_{\k-1}(K\cup Y;\Z)=0$ a.a.s.
\end{cor}

\begin{pf}
We have that the $(\k-1)$-skeleton of $\Y$ coincides with that of $K\cup Y$, and the set of $\k$-simplices of $\Y$ is contained in that of $K\cup Y$.
Thus $\tilde{H}_{\k-1}(K\cup Y;\Z)$ is a quotient of $\tilde{H}_{\k-1}(\Y;\Z)$ and so it follows from Lemma \pr{p8} that
$\tilde{H}_{\k-1}(K\cup Y;\Z)=0$ a.a.s. 
\end{pf}

\begin{lemma}\label{20}
If $\tilde{H}_{\k-1}(K\cup Y;\Z)=0$ then $Y\supseteq K$, so $K\cup Y = Y$.
\end{lemma}

\begin{pf}
Assume on the contrary that there exists an $(\k-1)$-simplex $\rho\not\in Y$. Let $\sigma$ be an $\k$-simplex such that $\rho \in \pa \sigma$. Then $\pa\sigma$ is an $(\k-1)$-cycle in $K\cup Y$ which cannot be a boundary in $K\cup Y$ since $\rho$ is not contained in any $\k$-simplex of $K\cup Y$.
(Note that if $\k-1=0$ then $\pa\sigma$ is indeed a \emph{reduced} cycle in $K\cup Y$.)
\end{pf}

This leads us to our two concluding propositions.

\begin{prop}\label{c3}
$\tilde{H}_{\k-1}(Y;\Z)=0$ a.a.s.
\end{prop}

\begin{pf}
By Corollary \pr{c2} we have $\tilde{H}_{\k-1}(K\cup Y;\Z)=0$ a.a.s. but then by 
Lemma \pr{20} we have $K\cup Y = Y$ so in fact $\tilde{H}_{\k-1}(Y;\Z)=0$ a.a.s.
\end{pf}

\begin{prop}\label{p7}
$Y$ contains the full $(\k-1)$-skeleton a.a.s. and so $\tilde{H}_k(Y;\Z)=0$ for all $k<\k-1$ a.a.s.
\end{prop}

\begin{pf}
By Corollary \pr{c2} we have $\tilde{H}_{\k-1}(K\cup Y;\Z)=0$ a.a.s. so by Lemma \pr{20} we have $Y\supseteq K$ a.a.s.
\end{pf}
For a different proof of Proposition \pr{p7} see Lemma 11.6 of \cite{fmn}.

Finally, Propositions \pr{p6},  \pr{dec}, \pr{bk}, \pr{c3},  \pr{p7} together constitute our desired Theorem \pr{p3}.
This is complemented by Proposition \pr{A} that covers the case $\beta<0$.

\end{document}